\documentclass[12pt,twoside,a4paper]{amsart}
\usepackage{amsfonts}
\usepackage{amssymb}
\usepackage{amsmath}
\usepackage{amscd}
\usepackage{color}
\usepackage[all]{xy}
\usepackage{mathrsfs}
\usepackage{xspace,yhmath,mathdots,txfonts}
\usepackage{yfonts}

\usepackage{hyperref}

\DeclareMathAlphabet{\mathpzc}{OT1}{pzc}{m}{it}

\begin{document}

\newcommand{\pa}{\partial}
\newcommand{\opa}{\overline\pa}
\newcommand{\ol}{\overline }

\numberwithin{equation}{section}

\newcommand\C{\mathbb{C}}
\newcommand\R{\mathbb{R}}
\newcommand\Z{\mathbb{Z}}
\newcommand\N{\mathbb{N}}
\newcommand\PP{\mathbb{P}}

\newcommand\crdim{{{CR}\,\mathrm{dim}}}
\newcommand\crcodim{{CR\,\mathrm{codim}}}

\newcommand\ind{{\mathrm{ind}}}
\newcommand\coind{{\mathrm{coind}}}
\newcommand\Lf{\mathbf{L}}
\newcommand\Ff{\mathbf{F}}
\newcommand\et{\mathfrak{e}}

\newcommand\Ob{\mathbf{O}}
\newcommand\SO{\mathbf{SO}}
\newcommand\ot{\mathfrak{so}}

\newcommand\csO{\mathcal{O}_{M,\mathrm{str}}}
\newcommand\csOs{\widetilde{\mathcal{O}}_{S}}
\newcommand\sOs{\hat{\mathcal{O}}_{S,\mathrm{str}}}
\newcommand\so{\widetilde{\mathcal{O}}}
\newcommand\sO{\widetilde{\mathcal{O}}_{M}}

\newcommand\Kf{\mathbf{K}}
\newcommand\kil{\mathpzc{q}}
\newcommand\Qf{\mathbf{Q}}
\newcommand\Sf{\mathpzc{Stab}}
\newcommand\sta{\mathfrak{stab}}
\newcommand\Tf{\mathbf{T}}
\newcommand\Tfp{\mathbf{T}^+}
\newcommand\Tfm{\mathbf{T}^-}
\newcommand\Gfu{\mathbf{G}_{u}}
\newcommand\gtu{\mathfrak{g}_u}
\newcommand\GL{\mathbf{GL}}
\newcommand\End{\mathpzc{End}\!}
\newcommand\Aut{\mathpzc{Aut}}
\newcommand\Wg{\mathpzc{W}}
\newcommand\wt{\mathfrak{w}}
\newcommand\ttu{\mathfrak{t}_u}
\newcommand\qt{\mathfrak{q}}
\newcommand\gt{\mathfrak{g}}
\newcommand\Cmp{\mathpzc{C}}
\newcommand\Rad{\mathpzc{R}}
\newcommand\hg{\mathfrak{h}}
\newcommand\Tq{\mathpzc{T}}

\newcommand\If{\mathbf{V}}
\newcommand\vt{\mathfrak{v}}

\newcommand\eg{\mathfrak{e}}
\newcommand\ig{\mathfrak{s}}
\newcommand\Vf{\mathbf{V}}

\newcommand\Zf{\mathpzc{Z}}
\newcommand\Zb{\mathbf{Z}}
\newcommand\Zbp{\mathbf{Z}^+}
\newcommand\Zbm{\mathbf{Z}^-}
 \newcommand\op{\mathrm{p}_0}
\newcommand\pr{\mathrm{pr}}
\newcommand\Qq{\mathpzc{Q}}

\newcommand\Wf{\mathbf{W}}
\newcommand\SL{\mathbf{SL}}
\newcommand\SU{\mathbf{SU}}
\newcommand\CP{\mathbb{CP}}
\newcommand\Gr{\mathpzc{Gr}}
\newcommand\Fi{\mathpzc{F}}

\newcommand\Xf{\mathrm{X}}
\newcommand\Yf{\mathrm{Y}}
\newcommand\Xg{\mathfrak{X}}
\newcommand\Ii{\mathpzc{I}}
\newcommand\Pro{\mathbb{P}}
\newcommand\Sb{\mathbf{S}}
\newcommand\Id{\mathrm{I}}
\newcommand\PC{\mathbb{PC}}
\newcommand\Ci{\mathcal{C}^\infty}
\newcommand\cC{\mathcal{C}}  
\newcommand\diag{\mathrm{diag}}
\newcommand\Bf{\mathbf{B}}
\newcommand\bt{\mathfrak{b}}
\newcommand\eq{\mathpzc{e}}

\newcommand\Nf{\mathbf{N}}
\newcommand\nt{\mathfrak{n}}
\newcommand\Uf{\mathbf{U}}
\newcommand\ut{\mathfrak{u}}
\newcommand\pg{\mathfrak{p}}
\newcommand\po{\mathfrak{p}_0}
\newcommand\Pf{\mathpzc{P}}
\newcommand\Pfo{\mathpzc{P}_0}
\newcommand\Mf{\mathpzc{M}}
\newcommand\Xff{\mathpzc{X}}
\newcommand\gl{\mathfrak{gl}}
\newcommand\trac{\mathrm{trace}}
\newcommand\tracr{\mathrm{trace}_{\mathbb{R}}}
\newcommand{\Ef}{\mathbf{E}}
\newcommand\Pj{\mathbb{P}}
\newcommand\slt{\mathfrak{sl}}
\newcommand\su{\mathfrak{su}}
\newcommand\Jd{\mathrm{J}}
\newcommand\Jf{\mathpzc{J}}
\newcommand\rt{\mathfrak{r}}
\newcommand\Pp{\mathfrak{P}}
\newcommand\td{\mathfrak{t}}
\newcommand\Cf{\mathbf{C}}
\newcommand\Ml{\mathpzc{M}}
\newcommand\cg{\mathfrak{c}}
\newcommand\Af{\mathbf{A}}
\newcommand\tg{\textgoth{t}}
\newcommand\HNR{$\mathrm{HNR}$}
\newcommand\kll{\mathpzc{k}}
\newcommand\Ifs{\mathbf{S}}
\newcommand\ift{\mathfrak{s}}
\newcommand\Tb{\mathbf{T}}
\newcommand\xt{\mathfrak{x}}
\newcommand\epi{\epsilonup}
\newcommand\aq{\mathpzc{a}}
\newcommand\bq{\mathpzc{b}}
\newcommand\rk{\mathrm{rank}}
\newcommand\cR{\mathcal{R}}
\newcommand\eps{\epsilon}
\newcommand\cM{\mathpzc{M}}
\newcommand\pct{\mathrm{p}}
\newcommand\qct{\mathrm{q}}
\newcommand\zct{\mathrm{z}}
\newcommand\esf{\textsf{E}}
\newcommand\iq{\mathpzc{i}}
\newcommand\pq{\mathpzc{p}}
\newcommand\qq{\mathpzc{q}}
\newcommand\cH{\mathcal{H}}
\newcommand\cF{\mathcal{F}}
\newcommand\Cq{\mathpzc{C}}
\newcommand\Cg{\mathfrak{C}}
\newcommand\Og{\mathfrak{O}}
\newcommand\xq{\mathpzc{x}}
\newcommand\zq{\mathpzc{z}}
\newcommand\cI{\mathcal{I}}
\newcommand\cN{\mathcal{N}}
\newcommand\vq{\mathpzc{v}}
\newcommand\wq{\mathpzc{w}}
\newcommand\uq{\mathpzc{u}}
\newcommand\sft{\mathpzc{stab}}
\newcommand\sfa{\textsf{a}}
\newcommand\sfd{\textsf{d}}
\newcommand\sfs{\textsf{s}}
\newcommand\Lq{\mathcal{L}}
\newcommand\Lb{\mathscr{L}}
\renewcommand\Re{\mathfrak{R}}
\newcommand\Rk{\mathcal{R}}
\newcommand\Xm{\mathfrak{X}}
\newcommand\gm{\textsf{g}}
\newcommand\hm{\textsf{h}}

\def\hksqrt{\mathpalette\DHLhksqrt}
   \def\DHLhksqrt#1#2{\setbox0=\hbox{$#1\sqrt{#2\,}$}\dimen0=\ht0
     \advance\dimen0-0.2\ht0
     \setbox2=\hbox{\vrule height\ht0 depth -\dimen0}
     {\box0\lower0.4pt\box2}}

\title{Aspects of the Levi form}

\author[J.~Brinkschulte]{Judith Brinkschulte}
\address{J.\ Brinkschulte:
Mathematisches Institut\\ Universit\"at Leipzig\\
Augustusplatz 10/11\\ 04109 Leipzig (Germany)}
\email{brinkschulte@math.uni-leipzig.de}

\author[C.D.~Hill]{C. Denson Hill}
\address{C.D.\ Hill:
Department of Mathematics\\ Stony Brook University
\\ Stony Brook, N.Y. 11794 (USA)}
\email{dhill@math.stonybrook.edu}

\author[J.~Leiterer]{J\"urgen Leiterer}
\address{J.\ Leiterer:
Institut f\"ur Mathematik\\ Humboldt-Universit\"at zu Berlin
Unter den Linden 6\\ 10099 Berlin
(Germany)}
\email{leiterer@mathematik.hu-berlin.de}

\author[M.~Nacinovich]{Mauro Nacinovich}
\address{M.\ Nacinovich:
Dipartimento di Matematica\\ II Universit\`a di Roma
``Tor Ver\-ga\-ta''\\ Via della Ricerca Scientifica\\ 00133 Roma
(Italy)}
\email{nacinovi@mat.uniroma2.it}

\maketitle

\vspace{0.5cm}

\begin{abstract}
 We discuss various analytical and geometrical aspects of the Levi form,
  which is associated with a CR manifold having any CR dimension and any
  CR codimension.
\end{abstract}

\vspace{0.5cm}

\section{Introduction}

The Levi form is a rather important geometric notion, which appears in
a fundamental way in several complex variables, complex differential
geometry, algebraic geometry, in certain aspects of partial differential
equations, and in particular in the theory of Cauchy-Riemann structures
on real manifolds (known for short as CR manifolds). Its role is to
measure certain second order effects, which are of natural interest in
those subjects. However it appears in various incarnations, and many
different authors have used it in different ways, each employing their
own peculiar notation, way of writing it, and their own understanding
of its geometric significance.\\

This has often led to confusion, in which even experts are unable to
easily decipher what some other researcher has written. This article
is an attempt to rectify the situation, and especially to establish
once and for all a good and consistent notation, which distinguishes
among the various incarnations of the Levi form. We hope that in the
future mathematicians will find it helpful and convenient to adopt
our conventions.\\

For the convenience of the reader, and to make everything understandable
to people not already familiar with the Levi form, we have begun with
a discussion of almost complex manifolds, passing to complex manifolds,
then to abstractly defined CR manifolds (or almost CR manifolds), and
finally to the situation of locally CR embeddable CR manifolds. And
along the way, we try to keep track of the various aspects of the Levi
form which naturally appear and explain the connections among them.\\

We have also stressed the various geometric meanings, of which there
are several. In the CR embedded case, we explain the connection of the
Levi form with the second fundamental form with respect to the induced
metric from the ambient space. Finally, for the important case of
homogeneous CR manifolds, we explain how the Levi form can be computed
from the point of view of Lie algebra, and illustrate it with a couple
of examples.

\section{Almost complex manifolds and the Nijenhuis tensor}

Let  $M$ be a smooth (real) manifold of dimension $2n$. An 
{\it almost complex structure} on $M$ is 
a smooth assignment of a complex structure on each fiber
of $TM,$ i.e. a fiber preserving smooth map $J\,{:}\,TM\,{\to}\,TM$
which is linear on the fibers and satisfies $J^2=-\mathrm{Id}.$
This map $J$ uniquely 
extends to a smooth vector bundle automorphism 
of the complexified tangent bundle 
$\C TM\ ({\,=}{\,\C}\otimes TM)$, 
which we 
denote by 
the same symbol
$J.$ 
The condition $J^2 = -\mathrm{Id}$ 
implies that $J$ has eigenvalues $i$ and $-i$, 
and we get a decomposition
\begin{equation}\label{15.7.18n}
\C TM = T^{1,0}M \oplus T^{0,1}M\end{equation}
into the corresponding eigenspaces
\begin{equation*}T^{1,0}M = \lbrace X- iJX \mid X\in TM 
\rbrace\end{equation*}
and
$$T^{0,1}M = \lbrace X+ iJX \mid X\in TM \rbrace.$$
We also note that 
\begin{equation}  \label{acx1}
T^{0,1}M = \ol{T^{1,0}M}\textcolor{yellow}{,}\textcolor{blue}{.}
\end{equation}

We call $T^{1,0}M$ the {\em $T^{1,0}$-bundle}, 
and $T^{0,1}M$ the  {\em $T^{0,1}$-bundle} 
of the almost complex structure $J$.
\\

To simplify notation, 
in the following 
we will use the same symbol
for  both vector fields and tangent vectors, whenever
we believe this will not cause confusion.\\

Equivalently, an almost complex structure can be defined by 
the assigning of
its $T^{1,0}$-bundle. 
Indeed, 
for  
each smooth complex subbundle $T^{1,0}M$ 
of $\Bbb C TM$ satisfying \eqref{15.7.18n} 
with 
$T^{0,1}M:=\overline{T^{1,0}M}$, 
there is a unique almost complex structure with this 
$T^{1,0}$-bundle. 
\par 
Since by \eqref{15.7.18n} the spaces $T_{x}^{1,0}M$ do not contain
nontrivial purely imaginary vectors, taking the
real part is an $\R$-linear isomorphism 
$\mathfrak{R}_{x}$ 
of $T_{x}^{1,0}M$
onto $T_{x}M$ and  
therefore we may define $J$ 
on $T_{x}M$ 
by requiring that
$X{-}iJX\in{T}_{x}^{1,0}M.$ The complex structure $J$
on $TM$ 
extends to an anti-involution $J^{\C}$ on $\C{T}M,$
which is multiplication by $i$ on  $T^{1,0}M$ and
$({-}i)$ on $T^{0,1}M$.
With the notation
\begin{equation*}
\mathfrak{R}:T^{1,0}M \to TM,\;\; \mathfrak{R}(X)=\tfrac{1}{2}
(X{+}\bar{X}),\;\forall X\in{T}^{1,0}M, 
\end{equation*}
we have \begin{equation*}
JX=\mathfrak{R}(\iq{\cdot}\mathfrak{R}^{-1}(X)),\;\;\forall
X\in{T}M.
\end{equation*}

In the same way, the $T^{0,1}$-bundle can also be used to define the almost complex structure, by giving a smooth complex subbundle $T^{0,1}M$ of $\Bbb C TM$ satisfying \eqref{15.7.18n} with $T^{1,0}M:=\overline{T^{0,1}M}$.

We refer the reader to \cite{B} for more details.\\

A complex manifold $M$ admits an almost complex structure,
which can be described by using its complex coordinate charts
to locally define its $T^{1,0}$-bundle to be the span of 
\begin{equation*}
\dfrac{\partial}{\partial{z}_j}\coloneqq \frac{1}{2}
\left(\frac{\partial}{\partial{x}_j}-
i\,\frac{\partial}{\partial{y}_j}\right),
\;\; j=1,\hdots,n,
\end{equation*} 
for any set of local holomorphic coordinates 
$z_1{=}x_1{+}i y_1,\hdots,z_n=x_n{+}i y_n,$
with  
underlying real coordinates $x_j,y_j$ ($j{=}1,\hdots,n$). 
By the holomorphic chain rule, this is a good definition
and we call this object \emph{the almost complex structure
defined by the complex structure of $M.$}\\

On the other hand, not every almost complex structure 
is defined by a complex structure.  
If $T^{1,0}M$ is the $T^{1,0}$-bundle 
of an almost complex structure defined by a complex structure, 
then the formal integrability condition 
\begin{equation}  \label{cx1}
\lbrack T^{0,1}M, T^{0,1}M \rbrack \subset T^{0,1}M
\end{equation}
is satisfied.
This means that if $X$ and $Y$ are smooth (real) vector fields, 
then the commutator of $X+iJX$ and $Y+iJY$ is of the same sort, 
i.e.
$$\lbrack X+iJX, Y+iJY\rbrack = Z+iJZ$$
for some smooth (real) vector field $Z$.
One easily 
checks 
that this formal integrability 
condition is equivalent to the vanishing 
of the Nijenhuis tensor, 
defined 
in \cite{NW} 
by
\begin{equation} \label{cx2}
\mathcal{N}(X,Y) = \lbrack X,Y\rbrack - \lbrack JX,JY\rbrack 
+ J\lbrack X, JY\rbrack + J\lbrack JX, Y\rbrack
\end{equation}
for $X,Y\in\Gamma(M,TM)$. 
The Newlander-Nirenberg theorem \cite{NN} asserts 
that formal integrability of an almost complex structure 
implies in fact integrability, 
i.e. there exists a complex structure 
which defines the almost complex structure. 
An improved version of this theorem requiring 
a minimal amount of smoothness 
of the almost complex structure 
can be found in \cite{HT}. 
\section{Embedded and abstract (almost) CR manifolds}\label{s3}
Let $M$ be a 
smooth 
manifold 
of 
real 
dimension $2n{+}k$. 
An 
{\it almost CR structure} 
$(HM,J)$ 
of type $(n,k)$ on $M$ 
consists of 
the data of 
a smooth
subbundle $HM\subseteq TM$ 
of fiber dimension $2n$ 
and a 
fiber preserving 
smooth  
vector bundle isomorphism 
$J: HM\longrightarrow HM$ satisfying $J^2 = -\mathrm{Id}$. An
\emph{almost CR manifold} of type $(n,k)$ 
is a 
smooth
manifold endowed with an almost CR structure of 
type $(n,k).$ 
The number $n$ is called the CR dimension  and $k$  the CR
 codimension of $M$. In the following we will avoid for simplicity, 
 whenever possible,
 to explicitly
 mention 
 regularity assumptions on the manifolds under consideration. 
 \par 
An almost complex structure on a real manifold $M$ 
of dimension $2n$ is therefore the same as an almost CR structure 
  of CR codimension zero.\par

As in the almost complex case, an almost CR structure of type $(n,k)$  on $M$ 
can be equivalently defined by the datum of a complex subbundle 
$T^{1,0}M$ of complex fiber dimension 
$n,$  of the complexified tangent bunde $\C{T}M,$ 
satisfying 
\begin{equation}\label{cr0}
 T^{1,0}M\cap T^{0,1}M = \{ 0\}, \;\;\;
 \text{ where $T^{0,1}M = \ol{T^{1,0}M},$}
\end{equation} 
(this is equivalent to \eqref{15.7.18n} if $k=0$). 
Since $T^{1,0}M$ does not contain purely imaginary vectors, 
the real parts of its
vectors form a real subbundle $HM,$ with fiber dimension $2n$ 
of $TM$ and
$J$ is defined on $HM$ by requiring 
that $X{-}iJX\,{\in}T^{1,0}M$ for all $X{\in}HM.$ 
We have $T^{1,0}M = \lbrace X- iJX \mid X\in HM \rbrace$,
$T^{0,1}M = \lbrace X+ iJX \mid X\in HM \rbrace$, 
and $T^{1,0}M\oplus T^{0,1}M=\C HM$, where $\C HM$ 
is the complexification of $HM$, i.e.  
the complex linear span of $HM$ in $\C TM$. 
\par
We call $T^{1,0}M$ the {\em $T^{1,0}$-bundle}, and $T^{0,1}M$ 
the  {\em $T^{0,1}$-bundle} of the almost CR structure $J$.
\par
An almost CR structure on $M$ is called a \emph{CR structure}  
if
the formal integrability condition
\begin{equation} \label{cr1}
\lbrack T^{0,1}M, T^{0,1}M\rbrack \subseteq T^{0,1}M, 
\end{equation}
is satisfied. Equation \eqref{cr1} means
that for   
$X,Y\in\Gamma(M,T^{1,0}M)$ 
their commutator $[X,Y]$ 
still 
belongs to 
$\Gamma(M,T^{1,0}M).$
Since
$$\lbrack X+iJX, Y+iJY\rbrack = \lbrack X,Y\rbrack 
- \lbrack JX,JY\rbrack + i\lbrace \lbrack X, JY\rbrack 
+ \lbrack JX, Y\rbrack \rbrace,$$
for $X, Y\in \Gamma(M,HM)$,
formal integrability can be reformulated 
in terms of real vector fields by the two conditions, 
the second one involving a 
Nijenhuis tensor on $HM$:
\begin{gather}\label{cr2}
 [X,Y]-[JX,JY]\in 
 \Gamma(M,HM)
 \qquad\qquad\qquad\\ 
 \label{cr3}
 \mathcal{N}(X,Y) {=} \lbrack X,Y\rbrack {-} \lbrack JX,JY\rbrack {+} 
 J(\lbrack X, JY\rbrack {+} \lbrack JX, Y\rbrack) {=} 0 
 \qquad\qquad\\
 \notag
\qquad\quad\qquad\qquad\qquad\qquad\qquad\;\forall X,Y\in 
\Gamma(M,HM).
\end{gather}
\par\smallskip
Prime examples of CR structures of CR dimension $n$ 
are provided by the real  
submanifolds $M$ 
of a complex manifold $\mathfrak{X}$ for which 
\begin{equation}   \label{dim}
 \dim_{\C}(\C{T}_{x}M\cap{T}^{1,0}_{x}\mathfrak{X})
 =n\;\;\;\text{for all $x\in{M}.$}
\end{equation}
Note that $\dim_{\R}M{\geq}2n$ and 
 $\dim_{\C}\mathfrak{X}\geq\dim_{\R}M-n.$
\par 
It is easy to check that 
\[
T^{1,0}M:=\big\{\C T_x M\cap T^{1,0}_x\mathfrak X\big\}_{x\in M}
\]
is a 
 smooth 
complex subbundle (of complex fiber dimension $n$) 
 of 
 $\C TM$  
satisfying  
\eqref{cr0}, \eqref{cr1}
and
therefore induces a CR structure $(HM, J_M)$ 
of type $(n,k)$
 on $M$, where, with ${\Re}Z:=\frac{Z+\overline Z}{2}$ 
 for $Z\in T^{1,0}\mathfrak X$,
\[
HM:={\Re}(T^{1,0}M)\quad\text{and}\quad J_M X
\coloneqq{\Re}(i{\Re}^{-1}X),\quad 
\forall X\in H_xM,\;x\in M.
\]
We say that
this  CR structure on $M$ 
is
\emph{induced} from the complex structure of $\mathfrak{X}$. 
If $J:T\mathfrak X\to T\mathfrak X$ 
is the almost complex structure of $\mathfrak X$ 
induced 
by its 
complex structure,  
i.e. 
\begin{equation*}
JX={\Re}(i{\Re}^{-1}X)\quad\text{for all}\quad X\in T\Xm,
\end{equation*}then 
$J_M$ is the restriction of $J$ to $HM$.
In this situation, we say that $M$  is \emph{CR embedded} in  
$\mathfrak{X}.$ 
\par
The embedding codimension $\ell$ of $M$ in $\mathfrak{X}$ is always
greater or equal to the CR codimension $k$ of $M$ and $\ell{-}k$ 
is even
(the complex dimension of $\Xm$ is $m{=}n{+}\tfrac{1}{2}(k{+}\ell)$).
\par
For each point $x_{0}$ of $M$ we can find an open neighbourhood $U$ 
of  $x_{0}$ in $\mathfrak{X}$ 
such that $M\cap{U}$ 
is given by 
\begin{equation} \label{em1}
M\cap{U}=\{x\in U\mid
 \rho_{1}(x)=0,\; \hdots,\;\rho_{\ell}(x)=0\}
\end{equation}
with 
smooth 
real valued 
functions $\rho_1,\ldots,\rho_{\ell}$  satisfying
\[
d\rho_1\wedge\ldots\wedge d\rho_{\ell}\not=0\quad\text{on}
\quad 
M\cap{V}.
\]
Then, for the $T^{1,0}$- and $T^{0,1}$-bundles, we have, 
for $x\in M\cap U$,
\begin{align}&\label{em2}
T_x^{1,0}M=\big\{Z\in T_x^{1,0}\Xm\,\big\vert\,\partial\rho_1(Z)
=\ldots=\partial \rho_{\ell}(Z)=0\big\},\\
&\label{em3}T_x^{0,1}M=\big\{Z\in T_x^{0,1}\Xm\,\big\vert\,
\overline\partial\rho_1(Z)=\ldots=\overline\partial 
\rho_{\ell}(Z)=0\big\}.
\end{align}
In particular, 
$M$ is generically embedded (near $x_{0}$) if and only if 
\begin{equation}\label{em4}
\pa\rho_1\wedge \ldots\wedge \pa\rho_{\ell} \neq 0  
\end{equation}
on a neighbourhood of $x_{0}$ in $V.$ 
\par\smallskip  
The notion of \textit{embedded} CR manifolds 
leads us to 
briefly discuss 
CR maps and CR embeddings 
in a more general setting:\par\noindent
If $(M, HM, J)$ is a CR manifold and $\mathfrak{X}$ 
is a complex manifold, then $M$ is {\em CR embeddable} 
into $\mathfrak{X}$ if 
one can find 
a 
smooth 
embedding $\Psi: M\longrightarrow \mathfrak{X}$ such that 
$\Psi(M)$ is CR embedded in $\Xm$ and 
its induced CR structure  
agrees
     with the pushforward by $\Psi$ 
     of the abstract CR structure 
     on 
     $M.$ 
     This means
     that $\Psi$ and $\Psi^{-1}$ are CR maps, 
 according to the definition below:    \\
Let $M_{1}$ and $M_{2}$ be CR manifolds. A 
smooth 
map $\Psi:M_{1}\,{\to}\,M_{2}$ 
is called CR if $\Psi_{*}(T^{1.0}M_1)\,{\subseteq}\,T^{1,0}M_{2}.$ 
An equivalent definition,
not involving the complexification of the differential, 
requires that 
\begin{equation}\label{crmap}
\Psi_{*}(HM_{1})\,{\subseteq}\,HM_{2}\;\;\;\text{and}\;\;\; 
\Psi_{*}{\circ}J_{1}\,{=}\,J_{2}{\circ}\Psi_{*}
\end{equation} 
 (where we set $J_{h}$ for the 
CR structure
on $M_{h},$ $h{=}1,2$).
If, moreover, $\Psi$ is a 
diffeomorphism and $\Psi_\ast(T^{1,0}M_1) = T^{1,0}M_2$, 
then $\Psi$ is called a CR {\em isomorphism}.\par
\medskip
A CR manifold is called {\em locally embeddable} 
if it admits local CR embeddings  into some 
complex 
Euclidean space.
One can prove (see, e.g., \cite[p.\ 118]{HTa1} or 
\cite[p.\ 187]{B})  that each locally embeddable 
$\Ci$-smooth  
CR manifold admits local CR embeddings 
whose images are generically embedded.\\

More details can be found in \cite{B} and \cite{BER}.\\
Real-analytic CR manifolds are always embeddable 
(see e.g. \cite{AH72, AF79}), but, 
in 
contrast to the complex situation (when $k=0$), 
the formal integrability condition (\ref{cr1}) is, 
for arbitrary 
$\Ci$-smooth $CR$ manifolds, 
not sufficient to provide 
local 
embedding into complex manifolds. 
The local CR embedding problem is in fact very difficult 
and will not be discussed in this survey.

\section{The Levi form of (almost) CR manifolds}

Let $M$ be a smooth real manifold with an almost 
CR structure of type $(n,k).$ 
The \emph{Levi form} of $M$ 
at $x$ 
is  
the  
Hermitian symmetric 
map
$$ \Lb^{1,0}_x : T^{1,0}_x M\times T^{1,0}_x M 
\longrightarrow \C T_x M/ \C H_xM $$
defined by $\Lb^{1,0}_x (Z_x,W_x) = \frac{1}{2i}\pi_x 
(\lbrack \ol Z,  W\rbrack_x)$. Here $\pi_x$
is the canonical projection $\pi_x: \C T_x M\rightarrow 
\C T_x M / \C H_xM$, and $Z,W$ 
are smooth sections of $T^{1,0}M$ extending $Z_x, W_x$. 
The class of $\lbrack \ol Z,  W\rbrack_x$ in $\C T_x M/ \C H_xM$ 
is in fact  
independent of the choice of the extensions. 
Indeed, let $L_1,\ldots,L_n$ be a smooth basis of $T^{1,0}M$ 
near $x.$  If $Z',W'$ are smooth sections of
$T^{1,0}M$ with $Z'_x=Z_x$ and $W'_x=W_x,$ then
$Z'=Z+{\sum}_{j}\alpha_{j}L_{j},$
$W'=W+{\sum}_{j}\beta_{j}L_{j},$
for smooth complex valued functions
$\alpha_j$, $\beta_j$  
vanishing at $0.$ Then 
the value at $x$ of  
\begin{align*}
 [\bar{Z}',W']-[\bar{Z},W]={\sum}_{j,k}\overline{\alpha}_j
  \beta_k[\overline{L_j},L_k]
 +{\sum}_{j}\bar{Z}(\beta_{j})L_{j}-
 {\sum}_{j}W(\bar{\alpha}_{j})\bar{L}_{j},
\end{align*}
is an element of $T^{1,0}_{x}M{+}T^{0,1}_{x}M\,{=}\,\C{H}_{x}M.$ 
The (fiberwise) linear projection 
\begin{equation*}
\Re:\C{TM}\ni X \longrightarrow \Re{X}=\tfrac{1}{2}(X{+}\bar{X})\in{T}M
\end{equation*}
of $\C{T}M$ onto $TM$ along $i{T}M$ 
induces by passing to the quotients a left inverse 
$\widehat{\Re}\,{:}\,\C TM/\C H \to TM/HM,$ 
of the inclusion 
$TM/HM\hookrightarrow{\C}TM/\C{HM}$ 
which factors 
$TM\hookrightarrow\C{TM}.$ 
\par
Let $x\in M$. Since ${\Re}$ is a real linear isomorphism 
from $T^{1,0}_x M$ onto $H_xM$, we have a
well-defined  real bilinear map $\Lb_x:H_xM\times H_xM\to T_xM/H_xM$ 
making the following diagram commute:
\begin{equation}
\label{l10}
\begin{CD}
T^{1,0}_x M \times T^{1,0}_x M @>{\Lb^{1,0}_x}>> \C T_x M/ \C H_xM \\
@VV{\Re}V     @VV{\widehat{\Re}}V \\
 H_xM\times H_xM @>{\Lb_x}>> T_x M /  H_xM.  \end{CD} 
 \end{equation}
 \par\smallskip
For 
$X,Y\in \Gamma(M,T^{1,0}M)$, we have, 
as $\overline{({\Re}\vert_{T^{1,0}_xM})^{-1}X_x}=X_x+iJX_x$ 
and $({\Re}\vert_{T^{1,0}_xM})^{-1}Y_x=Y_x-iJY_x$,
\[\begin{split}\Lb_x(X_x,Y_x)&:=\widehat{\Re}\Big(\frac{1}{2i}
\pi_x\big([X+iJX,Y-iJY]\big)\Big)\\
&=\widehat{\Re}\Big(\frac{1}{2i}\pi_x\Big([X,Y]+[JX,JY]
+i\big([JX,Y]-[X,JY]\big)\Big)\Big)\\
&=\frac{1}{2}\pi_x\big([JX,Y]-[X,JY]\big),
\end{split}\]
where, 
to obtain the last equality,  
we used 
the facts  
that $\pi_x$ is complex linear and 
that  
$\widehat{\Re}$ projects $\C T_xM/\C H_xM$ 
onto $T_xM/H_xM$ 
along $i(TM/HM).$ 
In particular, $\Lb_x$ is symmetric .
and $\Lb(JX_{x},JY_{x})=\Lb(X_{x},Y_{x}).$ 
\par\smallskip

If the almost CR structure satisfies \eqref{cr2}, then 
\begin{equation*}
 [JX,Y]\equiv \big( [J^{2}X,JY]=-[X,JY]\big)\;\;\;\mod{HM}
\end{equation*}
and hence 
\begin{equation}\label{17.7.18}
\Lb_x(X_x,Y_x)=\pi_x([JX,Y])=\pi_x([JY,X]).
\end{equation}
For the corresponding quadratic forms

\[
\Lq^{1,0}_x(Z_x):=\Lb^{1,0}_x(Z_x,Z_x)
\quad\text{and}\quad \Lq^{}_x(X_x):=\Lb^{}_x(X_x,X_x),
\] 
\eqref{l10} yields
the following commutative diagram: 
\begin{equation}\label{l1}
 \xymatrix@!0{T_{x}^{1,0}M \ar[drrrr]^{\Lq^{1,0}_{x}} \ar[dd]_{\Re}\\
 &&&& T_{x}M/H_{x}M \\
 H_{x}M \ar[urrrr]_{\Lq_{x}}
 }
\end{equation} 
Almost 
CR structures satisfying (\ref{cr2}) 
were already considered in \cite{T} 
and are called {\it partially integrable} (see also \cite{CS}).\par
\smallskip
The 
direction of  
$\Lq_x(X)$ does not change when $X$ is
  replaced by another  
  tangent vector in the real span 
  of 
  $\{X,JX\}$; indeed,
  if $Y = a X + b JX$, then 
	\begin{equation} \label{changeofbase}
	\Lq_x(Y) = (a^2 + b^2)\Lq_x(X).
	\end{equation} 
\par 
In many situations, it is useful to also consider scalar Levi forms 
(as 
considered e.g. 
in \cite{AFN} 
and
\cite{HN}). For this, we need 
to introduce  
the \emph{characteristic conormal bundle}, 
which is 
the annihilator 
$$
H^{0}{M}
= \lbrace \xi\in T^\ast M\mid \xi (X) = 0 \, 
\forall X\in \Gamma (M,HM)\rbrace.$$
of $HM$ in $T^\ast M.$
\par
 \smallskip
We 
define a family of scalar Levi forms parametrized 
by $\xi$ in the characteristic conormal bundle as follows:
Given $\xi\in 
H^{0}_{x}{M}$
 and $X_x,Y_x\in H_x M$, we choose smooth sections 
 $\tilde \xi$ of  
 $H^{0}{M}$  
 and $X,Y$ of $HM$ extending $\xi, Y_x,Y_x$.
By  
the invariant formula for the exterior derivative 
(see, e.g., \cite[p.\ 213]{Sp}), \[
d\tilde\xi(X,Y) = X( \tilde\xi(Y)) - Y( \tilde\xi(X)) 
- \tilde\xi(\lbrack X,Y\rbrack)
. 
\]
Since  $\tilde\xi(X)=\tilde\xi(Y)=0$,
this implies
$$d\tilde \xi (X_x,Y_x) = -\xi( \lbrack X,Y\rbrack ).$$
Hence both sides depend only on $\xi, X_x,Y_x$ 
and 
 we define 
\begin{equation}  \label{l2}
\Lq_x(\xi, X_x) = \xi (\lbrack J X, X\rbrack) = d\tilde\xi(X,JX),
\end{equation}
on $H_x M$. Note that $\Lq_x(\xi,\cdot)$ is  
Hermitian  
for the complex structure $J$ on $H_xM$: i.e.  
$\Lq_x(\xi,JX_x) = \Lq_x(\xi,X_x)$. \par
\smallskip
The image $\Rk_{x}$ of $\Lq_{x}$ is a real cone in $T_{x}M/H_{x}M.$ 
Its dual cone 
\begin{equation*}
 \Rk^{0}_{x}=\{\xi\in{H}^{0}_{x}M\mid \Lq_{x}(\xi,X)\geq{0},
 \;\forall X\in{H}_{x}M\}
\end{equation*}
 is a closed convex cone in $H^{0}_{x}M.$ 
\par 
It is also useful to consider 
the characteristic conormal sphere bundle 
$S^{0}M=(H^{0}M{\setminus}\{0\})/\R^{+}$  
whose typical fiber is the $(k{-}1)$-dimensional sphere.  
If $\varpi:{H}^{0}M{\setminus}\{0\}\to{S}^{0}M$ 
is the canonical projection,
 then $\varpi(\Rk^{0}{\setminus}\{0\})$ 
 is the total space of a bundle on $M$
 which has semialgebraic fibers
  and consist of a finite number of
  closed connected components if $M$ itself is semialgebraic. 

\par 

The scalar Levi forms $\Lq_x(\xi,X)$ enable us 
to microlocalize with respect to the Levi form $\Lq_x(X)$. 
A
geometric version of 
this 
microlocalization 
is discussed in 
\cite[\S{5}]{HN1}. 
\par  
When we think this would not cause confusion, we will drop the
subscript {\textquotedblleft{$x$}\textquotedblright} and simply
write $\Lq$ for $\Lq_{x}.$

\section{
Levi form of an embedded CR manifold}
Assume that $M$ is  
a $(2n{+}k)$-dimensional smooth real submanifold 
of an $m$-dimensional 
complex manifold $\Xm$   
(with $m{\geq}n{+}k$)
and  
that 
\eqref{dim} holds true, so that 
the complex structure of 
$\Xm$ 
induces 
on $M$ 
a  
$CR$ structure $(HM,J)$ of type $(n,k)$. 
 Our first goal
in this section 
  is 
 to express the Levi form in terms of defining 
functions for $M.$ 
To simplify notation
and better understand the invariant meaning  
of the construction, 
it is convenient to introduce the real operator 
\begin{equation*}d^{c}=i(\bar{\partial}-\partial),
\;\;\;\text{satisfying}\;\;\;
 \partial=\tfrac{1}{2}(d{+}i{\cdot}d^{c}),\;\; 
 \bar{\partial}=\tfrac{1}{2}(d{-}i{\cdot}d^{c}),\;\;
2 i{\cdot}\partial\bar{\partial}=d d^{c}.
\end{equation*}

By considering the local situation, we can assume that our
CR manifold $M,$  
of type $(n,k)$, 
is given by 
\begin{equation}\label{5.1}
M\cap U= \lbrace z\in U\mid \rho_1(z)= \ldots = \rho_{\ell}(z)=0 \rbrace,
\end{equation}
with $\ell{=}2m{-}2n{-}k{\geq}k,$ 
for a system 
$\rho_1,\ldots,\rho_{\ell}$  
of 
real $\Ci$ defining 
functions, 
defined on an open neighborhood $U$ of $x_o\in M$ in $\mathfrak{X}$ 
and 
satisfying
$d\rho_1\wedge\ldots\wedge d\rho_{\ell}\neq{0}$ on~$M$. 
By \eqref{em2}, the vectors 
$X$ in $HM$ 
are characterised by 
$d^{c}\rho_j(X)=0$  
for $1{\leq}j{\leq}\ell$ and therefore the pullbacks 
$\left.d^{c}\rhoup_{1}\right|_{TM},\hdots,
\left.d^{c}\rhoup_{\ell}\right|_{TM}$ of the $d^{c}$'s 
of the defining functions on $M$ 
span the characteristic bundle $H^{0}M.$ 
However they are not linearly independent if $\ell{>}k.$ 
\par\smallskip
For a section $X$ of $HM,$ 
with $X{=}\Re(Z){=}\tfrac{1}{2}(Z{+}\bar{Z})$ for $Z{\in}{T}^{1,0}M,$ 
we have by~\eqref{l2} 
\begin{align} \label{5.2a}
 \Lq(d^{c}\rho_{j},X)=dd^{c}\rho_{j}(X,JX)
 =2i\partial\bar{\partial}\rho_{j}(X,JX)=
\partial\bar{\partial}\rho_{j}(Z,\bar{Z}).
\end{align} 
Then \eqref{5.2a} yields 
\begin{equation}\label{1.10.18}
 \Lq(d^{c}\rho_{j},X)= {\sum}_{\mu,\nu 
 = 1}^{m}  \frac{\pa^2\rho_j(x)}{\pa z_\mu\pa\ol z_\nu}Z_\mu 
 \ol Z_\nu,
 \qquad 1\le j\le \ell.
\end{equation}
\par

A $\xi$ in $H^{0}_{x}M$ can be written as a linear combination $\xi={\sum}_{j=1}^{\ell}\xi_{j}d^{c}\rho_{j}(x)$
with real coefficients and  
correspondingly  
\begin{equation}  \label{7.10.18'}
\Lq_x(\xi, X ) = {\sum}_{j=1}^{\ell}{\sum}_{\mu,\nu = 1}^{m} 
\xi_j \frac{\pa^2\rho_j(x)}{\pa z_\mu\pa\ol z_\nu}Z_\mu \ol Z_\nu.
\end{equation}
Note that, 
when the embedding is generic ($\ell{=}k$)
the  
$d^c\rho_j$'s are a basis of $H^0M$ and the $\xi_j$
are uniquely determined by $\xi.$ \par 
 Formula (\ref{7.10.18'}) is actually 
 the classical definition of the scalar Levi forms 
 which 
was   
 used for a long time by various authors 
 (e.g. \cite{AFN}, \cite{ChSh}, \cite{HTa}, \cite{H}).
 \par\medskip
The Levi form yields important information on the geometry 
of the embedding $M\hookrightarrow\Xm.$ 
In fact its signature is related to mutual positions of $M$ 
and \textit{holomorphic balls}:
By a $p$-dimensional holomorphic ball we mean the image 
of a smooth embedding
$$F: \ol{\mathbb{B}}_p = \lbrace z\in \C^p\mid \vert z\vert
\leq 1\rbrace\longrightarrow \C^{n+k}$$ 
which 
is holomorphic on the open ball $\mathbb{B}_p$ 
($1$-dimensional holomorphic balls 
are usually called holomorphic discs).
 Let us first consider the case $k=1$ 
 (a piece of a real hypersurface in 
 $\Xm$). 
 In this case 
 $H^{0}_{x}M$ 
 is one-dimensional, so it is generated 
 by one characteristic conormal direction $\xi$. 
 Assume 
 that the scalar Levi form $\Lq_x(\xi,\cdot)$ 
 is nondegenerate and has $p$ positive and $q$ negative eigenvalues 
 (hence $p+q =n$). This means that we can find a 
 $p$-dimensional holomorphic ball lying on one side of $M$, 
 touching $M$  only with its center at $x$, tangent to the 
 $p$-dimensional positive eigenspace of $\Lq_x(\xi,\cdot)$
  and 
  a $q$-dimensional holomorphic ball lying on the other side 
  of $M$, also touching $M$ only with its center at $x$, 
 tangent to the $q$-dimensional negative eigenspace of 
 $\Lq_x(\xi,\cdot)$ (see 
e.g. 
 \cite[pp. 798-800]{AH}).\par 
 Vice versa, if $\Lq(\xi,\,\cdot\,)$ has for all nonzero
 $\xi$ in $H^0_xM$ at least $q$ negative eigenvalues,
 all holomorphic balls of dimension $d{\geq}q{+}\tfrac{1}{2}(\ell
 {+}k)$ centered at $x$ have an intersection with $M$  of
 positive dimension. 
\par\smallskip
Now let $k\geq 1$ be arbitrary, and let us assume that, 
for some $\xi\,{\in}\,H^{0}_{x}M,$ the 
scalar Levi form 
$\Lq_x(\xi,\cdot)$ 
has $p$ positive eigenvalues.  
Then 
one can find a holomorphic ball, of complex dimension 
$p{+}\tfrac{1}{2}(\ell{+}k){-}1,$
touching $M$ only with its center at $x$, where it is 
 tangent to 
 a complex 
 $p$-dimensional 
 complex linear subspace of $H_{x}M$ 
 on which $\Lq_{x}(\xi,\,\cdot\,)$ is positive definite.
 Indeed, assuming as we can that $\xi=d^{c}\rho(x)$ 
 for a smooth real valued function vanishing on $M,$
 for a large real $c$ the set 
 $\big\{\rho{+}c{\cdot}{\sum}_{j=1}^{\ell}\rho_{j}^{2}\,{=}\,0\big\}$ 
 defines near $x$ a smooth CR hypersurface
 having a Levi form with $p{+}\tfrac{1}{2}(\ell{+}k){-}1$ 
 positive eigenvalues.
\par\medskip
Next we show that, 
 after introducing a Hermitian metric on $\Xm,$
the Levi form can be rewritten
as a map with values in 
 the normal space 
 of $M$, as 
 was
 suggested in \cite{Her}. 
\par\smallskip
A Riemannian metric $\gm$ 
on $\Xm$ is  \emph{Hermitian} if the complex structure
$J$ is an isometry on $T\Xm.$ In particular
\begin{equation*}
\gm(X,JX)=0
\;\;\text{and}\,\;\gm(JX,JX)=\gm(X,X),\;\;\forall
X\in{T}\Xm.
\end{equation*}
This is equivalent to the fact that 
$\gm$ is the real part of a complex valued
Hermitian symmetric scalar product $\hm.$ The tensors $\gm$ and $\hm$ are related by
\begin{equation*}
\hm_x(X,Y)=\gm_x(X,Y)-i{\cdot}\gm_x(JX,Y), \;\;\forall x\in\Xm,\;
\forall
X,Y\in{T}_x\Xm.
\end{equation*}
\par
Let $NM$ be the \emph{normal bundle} of $M$ in $\Xm,$ whose fibre at
$x\,{\in}\,M$ is the space 
$(T_{x}M)^{\perp}\,{=}\,\{X_{x}\,{\in}\,T_{x}\Xg\,
{\mid}\,\gm(X_{x},Y_{x})=0,\;
\forall Y_x\,{\in}\,T_{x}M\}$  
and $\pi_{N}{:}T|_{M}\Xm\,{\to}\,NM$  
the orthogonal projection. 
For $X\,{\in}\,T|_{M}\Xm,$ we set $X^{NM}=\pi_N(X).$ 
\par 
Since $HM = TM\cap JTM$, we have   
\begin{equation*}
 HM\,{=}\,\{X\,{\in}\,TM\,{\mid}\,\pi_{N}(JX)\,{=}\,0\}
 \end{equation*} 
 and therefore we obtain a commutative
 diagram 
\begin{equation}\label{5.5a}
 \xymatrix@!0{TM \ar[dd] \ar[drrr]^{\pi_{N}\circ{J}} \\
 &&& NM \\
 TM/HM \ar[urrr]_{\quad\varpi}}
\end{equation}
where the vertical arrow is projection onto the quotient and $\varpi$ is injective.
Note that $\varpi$ 
is a linear isomorphism when the embedding is generic.
By composing $\Lq$ 
with $\varpi$ we obtain a Hermitian symmetric quadratic form
\begin{equation}\label{5.6a}
 \Lq^{NM}(X)= \varpi\circ\Lq(X)=
 (J\, [JX,X])^{NM},\;\;\text{for $X\,{\in}\,HM,$}
\end{equation}
on $HM$ 
taking value in the normal bundle. \par\smallskip 
If $M$ is given locally 
by 
a system of local defining functions  \eqref{5.1}, we can use \eqref{5.2a}  to obtain a
description of $\Lq^{NM}.$ Indeed, 
the \textit{gradient} $\nabla{f}$ 
of a real valued
smooth function $f$ is the real vector field characterised 
by $\gm(\nabla{f},X)\,{=}\,df(X)$ 
for all $X\,{\in}\,\Gamma(\Xm,T\Xm).$ Thus 
\begin{align*}
 \Lq(d^{c}\rho_{j},X)=d^{c}\rho_{j}([JX,X])=-d\rho_{j}(J\,[JX,X])=-\gm(\nabla\rho_{j},J[JX,X])
\end{align*}
and \eqref{1.10.18}  expresses  the scalar product of $\Lq^{NM}(X)$ with
the \textit{gradient} of the defining function $\rho_{j}.$ The matrix
$(\gm(\nabla\rho_{j},\nabla\rho_{h}))$ is symmetric and positive and, if $(A_{j,h})$ is its inverse,
we obtain 
\begin{equation}\label{5.6}
 \Lq^{NM}(X)=-{\sum}_{j,h} A_{j,h} \Lq(d^{c}\rho_{j},X)\,{\cdot}\,\nabla\rho_{h}.
\end{equation}
The defining functions $\rho_1,\ldots,\rho_{\ell}$ may be chosen 
in such a way that, at a point $x\,{\in}\,M,$  
their gradients form 
an orthonormal basis\footnote{This can be 
achieved by applying the Gram-Schmidt orthonormalising
process to  the $\ell$-tuple 
$(\nabla\rho_1(x),\ldots,\nabla\rho_\ell(x)).$ 
}
of $N_xM.$ Fix holomorphic coordinates $z_{1},\hdots,z_{m}$ on $\Xm$ centered at $x.$
Then, with $X_{x}{=}\Re{Z_{x}}$ for a $Z={\sum}_{\mu=1}^{m}Z_{\mu}(\partial/\partial{z}_{\mu})
\in{T}^{1,0}_{x}M,$
\eqref{5.6} yields
\begin{equation} \label{7.10.18}
{\Lq}^{NM}(X_x) = -{\sum}_{j=1}^{\ell}  \Big({\sum}_{\mu,\nu=1}^{m} 
\frac{\pa^2 \rho_j(x)}{\pa z_\mu \pa\ol z_\nu}Z_\mu\ol Z_\nu \Big) \nabla\rho_j(x).
\end{equation}\\
\par\medskip
There is a unique linear connection $\nabla$ for which both the metric tensor $\gm$ and 
the complex structure $J$ are parallel (see e.g. \cite{GH}). 
Let us assume 
moreover
that $\Xm$ is K\"ahler, so that the Hermitian 
and the Levi-Civita connections 
coincide. Denote by $\nabla$ the covariant derivative on $\Xm.$ 
If $X,Y$ are vector fields on $\Xm$ which are tangent to $M,$ then taking
the normal component 
\begin{equation}\label{5.8}
B_{x}(X,Y)\,{=}\,(\nabla_{X_{x}}Y)^{NM},\;\;\;\text{ for $x\,{\in}\,M,$}
\end{equation}
of the covariant derivative of $Y$ with respect to $X$ 
defines  
an $NM$-valued  
symmetric tensor on $TM,$ which is called the
\emph{second fundamental form} of the embedding of $M$ into $\Xm.$ 
Likewise, the tangential projection $\nabla^{M}_{X}Y$ is the covariant derivative
of the Levi-Civita connection on $M$ of $\gm|_{M}.$ 
The Levi form $\Lq^{NM}$ can be expressed by using the second fundamental form
(see \cite{EHS} or \cite{HTa}):
\begin{equation} \label{3rdformula}
\Lq^{NM}_x(X) = B_x(X,X) + B_x(JX,JX)
\end{equation}
Indeed, if $X$ is a real vector field on $\Xm$ whose restriction to $M$ takes values in $HM,$ 
we obtain, in view of \eqref{5.5a} amd \eqref{5.6a}, 
\begin{eqnarray*}
B_x(X,X) + B_x(JX,JX) =\big(\nabla_X X\big)^{NM} + \big( \nabla_{JX}JX\big)^{NM}\qquad\qquad\\
=\big( \nabla_X X + J \nabla_{JX} X\big)^{NM} 
=\big( J\nabla_{JX} X - J\nabla_X JX\big)^{NM}\qquad\qquad \\
=\big( J\lbrack JX, X\rbrack\big)^{NM}=\Lq^{NM}(X).
\end{eqnarray*}
Since the hermitian connection has no torsion on
the two-dimensional planes $\langle{X},{JX}\rangle$
(see e.g, \cite{KN}), 
the second fundamental form is meaningful on complex
tangent lines also in the non-K\"ahler situation,
\eqref{3rdformula} also applies in the more general case. 

\section{On the geometric meaning of the Levi form}
In the remainder of this section, we want to 
make some more comments on the geometrical
interpretation of the Levi form. 
Let us start by discussing  
its 
invariance 
under $CR$ diffeomorphisms: \par

Let $(M_1, T^{1,0}M_1)$ and $(M_2,T^{1,0}M_2)$ be CR manifolds 
and  $\Psi\,{:}\, M_1\,{\to}\, M_2$   
a CR map.  
By

\eqref{crmap}
its differential  
 $\Psi_\ast$ 
 factors to yield 
 a map 
$$\lbrack \Psi_\ast(x)\rbrack: T_x M_1/(HM_1)_x \longrightarrow T_x M_2/(HM_2)_x,$$
 between 
 the quotient spaces
and one can check that
$$\lbrack\Psi_\ast(x)\rbrack\circ\Lq_x^{M_1} = \Lq_{\Psi(x)}^{M_2}\circ \Psi_\ast(x).$$
In particular, 
the Levi form is invariant under CR diffeomorphisms. \par\smallskip

When $M$ is locally generically CR embedded into $\C^{n+k}$,  we can find lots 
of local biholomorphic mappings near a fixed point  $x\in M$. 
All of them induce local CR diffeomorphisms of a neighborhood of $x$ in $M$ 
onto a generic CR submanifold of codimension $k$ 
sitting in another Euclidean space of 
the same  
complex dimension $n{+}k$. 
Even though the image of $M$ under such a local $CR$ 
diffeomorphism may look quite different geometrically, its Levi form
is invariant. 
An illustration of this is the fact that a strictly pseudoconvex hypersurface 
may not be 
\textit{convex}, but it is locally CR
diffeomorphic to a hypersurface that is strictly convex, in the
elementary sense. Here by a strictly pseudoconvex hypersurface 
we
mean one whose Levi form is either positive 
or negative
definite.\par\smallskip

More generally, as follows from Sylvester's 
{\textquotedblleft{law of inertia}\textquotedblright} for Hermitian forms, 
it is the signature of the Levi form that is invariant;
by the {\textquotedblleft{signature}\textquotedblright} 
of any scalar Levi form, we mean its number of positive, 
zero and negative eigenvalues.
It is therefore impossible to change the signature of the Levi form 
by making any ambient biholomorphic change of 
coordinates.
Nonetheless a simple rescaling procedure permits 
to change the "size" of  a nonzero eigenvalue.\par\smallskip

Using formula (\ref{3rdformula}), one can relate the Levi form of M to the extrinsic
curvature of M, when M is CR embedded.
 To see this, we assume that $M$ is  embedded into $\C^{m}$ near some point $x\in M$. 
 We will consider 2-dimensional "ribbons" inside $M$ constructed as follows: 
 Let $\ell$ be a complex line in $\C^{m}$ which passes through $x$ and is tangent to $M$ at $x$. 
  The (local) image
of the exponential map $\exp: TM \longrightarrow M$  restricted to $\ell$
 is then  a little piece of a real 2-dimensional surface $\Sigma$
 passing through $x$; 
actually $\Sigma$ is obtained as the union of all geodesics 
starting 
from $x$ 
in all directions of $\ell$. Equivalently,
      $\Sigma$ is locally given by the 
      {\textquotedblleft{slice}\textquotedblright} $M \cap\mathrm{span}\{\ell,N_xM\}$.\par\smallskip

 Now assume that $\ell= \mathrm{span}\{ X,JX\}$, and $X\in H_xM$ is of length one.
Since $N_x M\subset N_x\Sigma$, we then obtain from \eqref{3rdformula}
\begin{eqnarray}   \label{fundsigma}
\Lq^{NM}(X) & = &  B(X,X) + B (JX,JX)\nonumber\\
& = & \big( B_\Sigma(X,X) + B_\Sigma(JX,JX)\big)^{NM},
\end{eqnarray}
where $B_\Sigma$ is the second fundamental form of $\Sigma$, 
whereas $B$ is the second fundamental form of $M$. From 
 (\ref{fundsigma}) we get the following geometric interpretation 
 of the Levi form of $M$: $\Lq^{NM}(X) $ is twice the projection 
 of the mean curvature vector of $\Sigma$ onto $N_xM$.\par\smallskip

We would like to emphasize that this mean curvature vector depends on $x$ 
and $\ell$, but not on a particular choice of $X$. Indeed, if $\ell$ is spanned 
by $Y$ and $JY$, where $Y = aX{+}bJX$ and $a^2 {+} b^2 =1$, 
then we obtain the same mean curvature vector. 
In fact $\Lq^{NM}(Y) = \Lq^{NM}(X)$, as follows from (\ref{changeofbase}).\par\smallskip

  In contrast to the Gauss curvature of a surface (which is a metric invariant 
  but may change under conformal maps), 
  the mean curvature is known to be a conformal invariant: 
  it is possible to turn a strictly convex surface 
  looking like a skull-cap  into a saddle while keeping the mean curvature, say positive.\par\smallskip

 {\bf An Example: tube CR manifolds.} 
Let $\Gamma$ be a smooth real manifold in $\R^{n+k}$, of real dimension $n$ and real codimension $k$. Then $M= \Gamma\times \R^{n+k}$ is a CR manifold of type $(n,k)$, generically CR embedded into $\C^{n+k} = \R^{n+k} + i\R^{n+k}$. Indeed, if $pr: M\longrightarrow \Gamma$ denotes the canonical projection, then $T^{1,0}_xM = \C\otimes T_{pr(x)}\Gamma$. Following \cite{HK}, $M$ is called a tube CR manifold over the base $\Gamma$. \\

From (\ref{3rdformula}) we can deduce that the Levi form of $M$ at $x$ is the second fundamental form of $\Gamma$ at $pr(x)$: If we choose any real tangent vector $X$ to $\Gamma$ at $pr(x)$, then $X$ is also tangent to $M$ at $x$, and $JX = iX$ points along the flat fiber of $M$. Hence (\ref{3rdformula}) implies
\begin{eqnarray*}
\Lq^{NM}(X)&  = & B_x(X,X) + B_x (JX,JX) = B_x(X,X)\\
& = & \big( B^\Gamma_{pr(x)} (X,X) \big)^{NM} = \big( B^\Gamma_{pr(x)} (X,X) \big)^{N\Gamma}=
B^{\Gamma}_{pr(x)}(X,X).
\end{eqnarray*}
Therefore, assuming $X$ has length one, $\Lq^{NM}(X) =B^\Gamma_{pr(x)}(X,X)$ is the normal curvature vector to $\Gamma$ at $pr(x)$ in the direction $X$, and its length is the curvature of any curve in $\Gamma$ passing through $pr(x)$ having tangent vector $X$ at $pr(x)$.\par\smallskip

A similar observation applies to so-called Reinhardt CR manifolds.
 They are of the form $M = \exp{(\Gamma + i\R^{n+k})}$, 
 where $\exp$ denotes the exponential map 
 that sends $z= (z_1,\ldots, z_{n+k})$ to 
 $\exp z= (e^{z_1},\ldots, e^{z_{n+k}})$. 
 Reinhardt CR manifolds are thus invariant 
 under multiplication of each coordinate by a different 
 factor of modulus one, whereas tube CR manifolds are invariant 
 under translations in the pure imaginary direction. 

\section{The homogeneous case}

Using a notation which is customary while dealing with
homogeneous spaces (see e.g. \cite{Hel78})
we
add the subscript {\textquotedblleft{$0$}\textquotedblright}
to indicate \textit{real} objects, while \textit{complex} object
will be left
unsubscripted.\par  For a smooth bundle $\pr\,{:}\,E\to{M}$
we denote by $\Gamma(M,E)$ the space
$\{f\in\Ci(M,E)\mid \pr{\circ}f=\mathrm{id}_{M}\}$ 
of smooth sections of $E$
over $M.$ 
\subsection{Homogeneous CR manifolds and CR algebras} 
\paragraph{\textsc{Homogeneous spaces}}
Let $M_0$ be a smooth real
manifold and $\mathbf{G}_0$ a real Lie group, with identity $\eq.$
A differentiable 
action of $\mathbf{G}_0$ on $M_0$ is a $\Ci$ smooth map 
\begin{equation}\label{action}
\mathbf{G}_0\,{\times}\,M_0\ni (\sfa,\pct) \to \sfa{\cdot}\pct\in{M}_0,
\;\;
\text{with}\;\; \begin{cases} \eq{\cdot}\pct=\pct, 
\\
\sfa_1{\cdot}(\sfa_2{\cdot}\pct)=(\sfa_1{\cdot}\sfa_2){\cdot}\pct, 
\\ \forall\sfa_1,\sfa_2\in\mathbf{G}_0,\;
\forall \pct\in{M}_0.
\end{cases}
\end{equation}
\par
To each element $X$ of the Lie algebra $\gt_0$ of $\mathbf{G}_0$
we associate the smooth vector field $X^*\in\Gamma(M_0,TM)$ 
generating the flux $(t,\pct)\,{\to}\,\exp(t{\cdot}X){\cdot}\pct$
on $M_0.$ The space of $X^*$ for $X$ in $\gt_0$ is a 
real Lie subalgebra of $\Gamma(M_0,TM)$ for the commutation of
vector fields. The correspondence $X\to{X}^*$ is indeed
an antihomomorphism of Lie algebras, i.e. is $\R$-linear and
\begin{equation*}
[X^*,Y^{*}]=-[X,Y]^*,\;\;\forall X,Y\in\gt_0.
\end{equation*} 
\par 
 For each 
$\sfa\,{\in}\,\mathbf{G}_0,$ the map $L_\sfa(\pct)=\sfa\,{\cdot}\,\pct$ 
is a diffeomorphism of $M_0,$ 
yielding by differentiation a 
diffeomorphism $\sfd{L}_\sfa:TM_0\to{T}M_0.$
We will write for simplicity 
$\sfa\,{\cdot}\,\vq$ instead of $\sfd{L}_\sfa(\vq)$ if 
$\sfa\,{\in}\,\mathbf{G}_0$ and $\vq{\in}TM_0,$ noting that  
also the map
\begin{equation*}
\mathbf{G}_0\times{T}M_0\ni (\sfa,\vq)\to \sfa{\cdot}\vq\in{T}M_0
\end{equation*}
is a smooth action of $\mathbf{G}_0$ on $TM_0.$ 
\par
We say that $M_0$ is $\mathbf{G}_0$-homogeneous when the action
of $\mathbf{G}_0$ is transitive, i.e. when 
$\mathbf{G}_0\,{\cdot}\,\pct=M_0$ for some, and hence for all,
$\pct$ in $M_0.$ 
When this is the case, $T_\pct{M}_0=\{X^*_\pct\,{|}\,X{\in}\gt_0\}$
for all points $\pct$ of $M_0$ and therefore
one can use vector fields $X^*$ to get local frames
near any $\pct$ in $M_0.$\par\smallskip
The stabilizer 
 $\Sf_{\mathbf{G}_0}(\op){=}\{x{\in}\mathbf{G}_0\,{\mid}\, x{\cdot}\op{=}\op\}$ 
 of any point $\op$ of $M_0$ is a closed Lie subgroup
 of $\mathbf{G}_0.$ Denote by  
$\sft_{\gt_0}(\op)$ its Lie algebra.\par 
If the action  \eqref{action} is transitive, 
then $M_0$ has a unique real analytic structure such that 
for each point $\op$ the 
map 
\begin{gather}\label{subm}
\piup:\mathbf{G}_0\ni{\sfa}\to\sfa{\cdot}\op\in{M}_0,   
\intertext{is a real analytic submersion and  
the injective quotient 
yields a real analytic 
diffeomorphism}
\hat{\piup}:\mathbf{G}_0/\Sf_{\mathbf{G}_0(\op)}{\to}{M}_0\intertext{ and, by differentiating, a
linear isomorphism}  d\hat{\piup}_\eq:\gt_0/\sft_{\gt_0(\op)}\to
T_{\op}M_0.
\end{gather} 
\par\medskip

\paragraph{\textsc{Homogeneous CR manifolds}}
Let $M_0$ is a CR manifold with structure bundle $HM$
and partial complex structure $J.$ 
\par 
We say that a differentiable
action \eqref{action} is CR if the diffeomorfphisms  
$\pct\,{\to}\sfa{\cdot}
\pct$ of $M_0$ are CR for all $\sfa$ in $\mathbf{G}_0.$ 
This means that \begin{equation*}
\sfa\,{\cdot}\,HM=HM \;\;\text{and}\;\; \sfa\,{\circ}\,J\,{=}\,
J\,{\circ}\,\sfa.
\end{equation*}
This can be rephrased by using the complexification 
of the
action of $\mathbf{G}_0$ on $T{M}_0,$ by requiring that either one
holds true \begin{equation*}
\sfa{\cdot}T^{1,0}M_0=T^{1,0}M_0,\;\;\text{or}\;\;
\sfa{\cdot}T^{0,1}M_0=T^{0,1}M_0,\;\; \forall\sfa\in\mathbf{G}_0.
\end{equation*}
\par 
We say that $M_0$ is a 
\emph{$\mathbf{G}_0$-homogeneous CR manifold}
if the  
action \eqref{action} of the Lie group $\mathbf{G}_0$ on
$M_0$ is transitive and CR.
\par

Fix a base point $\op$ in $M_0$ 
and consider the smooth submersion \eqref{subm}.
The pullback 
\begin{equation*}
\Qq{=}{\piup}^*(\Gamma(M_0,T^{0,1}M_0)){=}
\{\zetaup \in \Gamma(\mathbf{G}_0,\C{T}\mathbf{G}_0\mid 
\sfd{\piup}(\zetaup_\sfa){\in}{T}^{0,1}_{\sfa{\cdot}{\op}},
\;\forall{\sfa}{\in}\mathbf{G}_0\}\end{equation*}
of $\Gamma(M_0,T^{0,1}M_0)$ 
 to $\mathbf{G}_0$ is a formally integrable distribution 
 of  complex  vector fields on $\mathbf{G}_0,$
 which is invariant by left translations.\par 
Let $\gt$ be the complex Lie algebra obtained by
complexifying  $\gt_0$ and set 
\begin{equation}\label{qu}
 \qt=\{Z\in\gt\mid \sfd\piup_e(Z)\in{T}^{0,1}_{\op}M_0\}.
\end{equation}
For the following statements we refer to \cite{AMN06,AMN07,MN05}.
\newtheorem{prop}{Proposition}[section]
\begin{prop} \label{propq} \begin{enumerate}
\item The subspace  
 $\qt$ defined in \eqref{qu}
 is a complex Lie subalgebra of $\gt.$ 
 \item $\qt\,{\cap}\,\bar{\qt}$ equals the complexification
 of $\sft_{\!\gt_0}(\op).$ 
\item A complex vector field  $\zetaup$ on $\mathbf{G}_0$ belongs to
 $\Qq$ iff $\sfa^{{-}1}_*\zetaup_\sfa\,{\in}\,\qt$ 
 for all $\sfa\,{\in}\,\mathbf{G}_0.$ 
 \item 
 Vice versa, any complex Lie subalgebra $\qt$ of $\gt$
 with $\qt\,{\cap}\,\bar{\qt}{=}
 \C{\otimes}_{\R}\sft_{\gt_0}(\op),$ 
defines on a $\mathbf{G}_0$-homogeneous space
$M_0$ a $\mathbf{G}_0$-ho\-mo\-ge\-neous CR structure
by setting 
\begin{equation*}\vspace{-19pt}
T^{0,1}_{\sfa{\cdot}\op}{M}_0=\sfa\,{\cdot}\,
\sfd\piup_\eq(\qt),\;\; \forall\sfa{\in}\mathbf{G}_0.
\end{equation*}
\qed
 \end{enumerate}
\end{prop}

A \textit{$\mathbf{G}_0$-homogeneous CR structure} on $M_0$
is determined by the datum of the antiholomorphic tangent
space at a point. Thus, 
by Proposition~\ref{propq}, 
having fixed a point $\op,$ 
the $\mathbf{G}_0$-homogeneous CR structures 
on $M_0$ are
parameterized by the 
\textit{complex} Lie subalgebras $\qt$ 
of the complexification $\gt$ of $\gt_0$ 
for which  
 \begin{equation}\qt\,{\cap}\,\gt_0
 =\sft_{\!\gt_0}(\op).\end{equation} 
A different choice of the base point $\op$ changes 
$\Sf_{\mathbf{G}_0}(\op),$ 
 $\sft_{\gt_0}(\op)$ and $\qt$ by $\mathbf{G}_0$-conjugation. 
 \par
These remarks motivate the definiton in \cite{MN05}
of  \emph{CR-algebras} as
pairs $(\gt_0,\qt)$ consisting of 
a \textit{real}
Lie algebra $\gt_0$ and a complex Lie subalgebra $\qt$
of   $\gt\,{=}\,\C\,{\otimes}_{\R}\gt_0.$ \par 
Indeed this notion encodes the possible realizations of $M_{0}$
as a smooth CR submanifold of a homogeneous complex manifold. 
Assume that: 
\begin{itemize}
\item $\mathbf{G}_0$ admits a complexification $\mathbf{G},$ 
\item there is a closed complex Lie subgroup
 $\Qf$ of $\mathbf{G}$ with Lie algebra $\qt$ and
 $\Qf\,{\cap}\,\mathbf{G}_0=\Sf_{\mathbf{G}_0}(\op).$
\end{itemize}
Then the inclusions $\mathbf{G}_0\hookrightarrow\mathbf{G}$ and
$\Sf_{\mathbf{G}_0}(\op)\hookrightarrow\Qf$ yield a 
smooth immersion of 
$M_0$ 
into the complex homogeneous space $M=\mathbf{G}/\Qf$ and we have
a commutative diagram 
\begin{equation*}
\begin{CD} \mathbf{G}_0 @>>> \mathbf{G} \\
@VVV @VVV \\
M_0 @>>> M.
\end{CD}
\end{equation*}
If we are only interested in the \textit{local} setting, we can always assume that 
the Lie group $\mathbf{G}_{0}$ is 
linear, thus admitting a complexification $\mathbf{G},$ and that $\Sf_{\mathbf{G}_{0}}(\op)$ 
is connected. 
 However, 
the analytic Lie subgroup $\Qf$ generated by $\qt$ may not be closed,
(\textit{virtual} in the sense of 
\cite{GOV93}). The quotient $\mathbf{G}/\Qf$ can be thought in this case
as a \textit{germ} of smooth complex manifold providing a local embedding of $M_{0}$
near $\op.$

\subsection{The Levi form} 
For a $\mathbf{G}_0$-homogeneous CR manifold $M_0$ one can describe
the Levi form at a point $\op$ by utilizing the attached 
CR algebra $(\gt_0,\qt)$ 
and Lie brackets on the complexification $\gt$ of $\gt_{0}.$
In fact we will consider the sesquilinear form $\Lq^{1,0}_{\op}$ on $T^{1,0}_{\op}M_{0}$
and the corresponding \textit{scalar} Levi forms, indexed by $H^{0}_{\op}M_{0}.$ 
With the submersion $\piup$
of \eqref{subm}
we have indeed $T^{1,0}_{\op}{M}_0=\sfd\piup_\eq(\bar{\qt}),$
where conjugation in $\gt$ is taken with respect to the real
form $\gt_0.$ 
The vector fields of the 
complex distribution $\Qq$ generated by the left-invariant vector fields
$Z^{*}$ with $Z\,{\in}\,\qt$ are $\piup$-related to the elements of $\Gamma(M_{0},T^{0,1}M_{0}).$ 
Hence we have (cf. e.g. \cite[Ch.I, Prop.3.3]{Hel78}) 
\begin{equation*}
 d\piup_{\op}([Z_{1},\bar{Z}_{2}])=[X_{1},\bar{X}_{2}]_{\op}
\end{equation*}
if $Z_{1},Z_{2}\in\Qq,$ $X_{1},X_{2}\in\Gamma(M_{0},T^{1,0}M_{0})$ and
$\piup_{*}(Z_{i})=X_{i},$ $i{=}1,2.$  
Hence
one obtains the following. 
\begin{prop} Let $\pr:\gt\to\gt/(\qt+\bar{\qt})$ 
be the projection 
onto the quotient  and 
$\varpi$ the isomorphism of $\gt/(\qt+\bar{\qt})$
with $\C{T}_{\op}M_0/\C{H}_{\op}M_0$ 
defined by the  isomorphism 
$\sfd\piup_\eq:\gt\to\C{T}_{\op}M_0,$
with $\sfd\piup_{\eq}(\qt+\bar{\qt})=\C{H}_{\op}M_0.$ \par 
Then the Levi form of $M_0$ at $\op$ is characterized
by the commutative diagram
\begin{equation}
\xymatrix{\qt\times\qt \ar[rr]^{(Z_1,Z_2)\to[Z_1,\bar{Z}_2]} 
\ar[d]^{(Z_1,Z_2)\to(\sfd\piup_{\eq}(\bar{Z}_1,
\sfd\piup_{\eq}(\bar{Z}_1))} 
&& \gt \ar[rr]^{\pr} &&\gt/(\qt+\bar{\qt}) \ar[d]^{\varpi}\\
T^{1,0}_{\op}M_0\times
T^{1,0}_{\op}M_0 \ar[rrrr]^{\Lq^{1,0}_{\op}} &&&& 
\C{T}_{\op}M_0/\C{H}_{\op}M_0
}
\end{equation}
where the top left are the Lie brackets 
in $\gt.$ \qed
\end{prop} 
We illustrate by some examples the actual computation. We need to
parameterise the quotient $T_{\op}M_{0}$ by taking a linear complement
of $(\qt\,{+}\,\bar{\qt})\,\cap\,\gt_{0}$ in $\gt_{0}.$ Then we will write the
Levi form $\Lq^{1,0}_{\op}$ as an $n{\times}n$ Hermitian symmetric matrix
(with $n$ equal to the CR dimension) depending upon $k$ real parameters
(with $k$ equal to the CR codimension). The Hermitian symmetric forms
obtained by fixing the real values of the parameters are the \textit{scalar}
Levi forms.
\theoremstyle{definition}
\newtheorem{exam}[prop]{Example}
\begin{exam}
 We consider the flag manifold $M$ consisting of the pairs $\ell_1{\subset}\ell_4{\subset}\C^6,$ 
 which is a  
compact complex manifold of complex dimension $11,$
homogeneous space for the action of $\SL_6(\C).$ Fixing on $M$ the base point
$\op{=}(\langle\eq_1\rangle,\langle
\eq_1,\eq_2,\eq_3,\eq_4\rangle),$ its stabiliser $\Qf$ in $\SL_6(\C)$ is the 
complex parabolic subgroup
with Lie algebra 
\begin{equation*}
\qt{=}\left\{\left.  
\begin{pmatrix}
 z_{1,1} & z_{1,2} & z_{1,3} & z_{1,4} & z_{1,5} & z_{1,6} \\
  0 & z_{2,2} & z_{2,3} & z_{2,4} & z_{2,5} & z_{2,6} \\
0 & z_{3,2} & z_{3,3} & z_{3,4} & z_{3,5} & z_{3,6} \\
 0 & z_{4,2} & z_{4,3} & z_{4,4} & z_{4,5} & z_{4,6} \\
0 & 0 & 0 & 0 & z_{5,5} & z_{5,6} \\
0 & 0 & 0 & 0 & z_{6,5} & z_{6,6} 
\end{pmatrix}\right| \text{with ${\sum}_{i{=}1}^6z_{i,i}=0$ }\right\}.
\end{equation*}
Let us consider the Hermitian symmetric matrix
\begin{equation*}
 B= 
\begin{pmatrix}
 0 & 0 & 0 & 0 & 0 & 1\\
  0 & 0 & 0 & 0 & 1 & 0\\
   0 & 0 & 1 & 0 & 0 & 0\\
    0 & 0 & 0 & 1 & 0 & 0\\
     0 & 1 & 0 & 0 & 0 & 0\\
      1 & 0 & 0 & 0 & 0 & 0
\end{pmatrix},
\end{equation*}
of signature $(4,2),$ to define 
$\SU(2,4){=}\{x{\in}\SL_6(\C) \,{\mid}\, x^*Bx=B\}.$ 
Then $\langle\eq_1\rangle$ is an isotropic line and 
$\langle\eq_1,\eq_2,\eq_3,\eq_4\rangle$ a $4$ plane 
on which the restriction of $B$ has minimal rank.
Hence the $\SU(2,4)$-orbit $M_0$ of $\op$ in $M$ is the minimal one. 
\par

The Lie algebra of $\SU(2,4)$ is characterized by 
\begin{equation*}
 \su(2,4)=\{X{\in}\slt_6(\C)\mid X^*B+BX{=}0\}.
\end{equation*}
Taking into account that $B^2=\Id_6,$ it is the subalgebra of fixed points 
of the conjugation $Z{\to}{-}BX^*B.$ 
Using this conjugation we obtain 
\begin{equation*}
\bar{\mathfrak{q}}{=}\left\{\left.  
\begin{pmatrix}
 z_{1,1} & z_{1,2} & z_{1,3} & z_{1,4} & z_{1,5} & z_{1,6} \\
  z_{1,2} & z_{2,2} & z_{2,3} & z_{2,4} & z_{2,5} & z_{2,6} \\
0 & 0 & z_{3,3} & z_{3,4} & z_{3,5} & z_{3,6} \\
 0 & 0 & z_{4,3} & z_{4,4} & z_{4,5} & z_{4,6} \\
0 & 0 & z_{5,3} & z_{5,4} & z_{5,5} & z_{5,6} \\
0 & 0 & 0 & 0 & 0 & z_{6,6} 
\end{pmatrix}\right| \text{with ${\sum}_{i{=}1}^6z_{i,i}=0$ }\right\}.
\end{equation*}
To compute the CR type and the Levi form of $M_0$ 
it is convenient to start from a linear complement
of $\qt$ in $\slt_6(\C)$ and construct the tangent 
to $M_0$ at $\op$ by adding, for each coefficient, its
$\su(2,4)$-conjugate.\par
We obtain a representation of the tangent space to $M_0$ 
at $\op$ of the form
\begin{equation*} \mathfrak{m}_0 \; : \; 
\begin{pmatrix}
 0 & 0 & 0 & 0 & 0 & 0\\
 \zq_1 & 0 & 0 & 0 & 0 & 0\\
 \wq_1 & {-}\bar{\zq}_2& 0 & 0 & 0 & 0 \\
 \wq_2 & {-}\bar{\zq}_3 & 0 & 0 & 0 & 0 \\
 \wq_3 & \iq{t}_2 & \zq_2 & \zq_3& 0 & 0 \\
 \iq{t}_1 & {-}\bar{\wq}_3 & {-}\bar{\wq}_1 & {-}\bar{\wq}_2 & {-}\bar{\zq}_1 & 0
\end{pmatrix}.
\end{equation*}
The $\zq_i,$ which are the terms in the complement of $\qt$ which by conjugation are sent into a
${-}\bar{\zq}_i$ lying in $\qt,$ correspond to the analytic tangent and should be therefore considered
{\textquotedblleft{complex}\textquotedblright} coordinates. The remaining terms, whose conjugates stay
in the complement of $\qt,$ are the {\textquotedblleft{real}\textquotedblright} coordinates, 
whose entries will parametrise the  Levi forms.  
The number of $z_i$'s is the CR dimension, twice the number of $\wq_i$'s plus the number or
$t_i$'s is the CR codimension:
in this example, $M_0$ is of type $(3,8)$ (we have $\zq_{1},\zq_{2},\zq_{3}$ and
$\wq_{1},\wq_{2},\wq_{3},t_{1},t_{2}$).
The Levi form is $3{\times}3$ and 
depends on $3$ complex (related to $\wq_1,\wq_2,\wq_{3}$)
and two real (corresponding to $\iq{t}_{1},\iq{t_2}$)
parameters. It can be represented by the matrix 
\begin{equation*} \tag{$*$}
\begin{pmatrix}
 0 & \wq_{1} & \wq_{2}\\
 \bar{\wq}_{1}& t_{2} & 0\\
 \bar{\wq}_{2} & 0 & t_{2}
\end{pmatrix}.
\end{equation*}
Let us explain the way it is computed and thus its meaning.
The $T^{1,0}_{\op}M_0$ part is represented in
$\mathfrak{m}{=}\C{\otimes}\mathfrak{m}_0$ by the matrices 
\begin{equation*}\tag{$**$}
 Z {=} 
\begin{pmatrix}
 0 & 0 & 0 & 0 & 0 & 0\\
 z_1 & 0 & 0 & 0 & 0 & 0\\
 0 & 0 & 0 & 0 & 0 & 0\\
 0 & 0 & 0 & 0 & 0 & 0\\
 0 & 0 & z_2 & z_3 & 0 & 0\\
 0 & 0 & 0 & 0 & 0 & 0
\end{pmatrix},
\end{equation*}
the $T^{0,1}M_0$ part by the matrices 
\begin{equation*}
 \bar{Z}{=}\begin{pmatrix}
 0 & 0 & 0 & 0 & 0 & 0\\
 0 & {-}\bar{z}_2 & 0 & 0 & 0 & 0\\
 0 & {-}\bar{z}_3 & 0 & 0 & 0 & 0\\
 0 & 0 & 0 & 0 & 0 & 0\\
 0 & 0 & 0 & 0 & 0 & 0\\
 0 & 0 & 0 & 0 & {-}\bar{z}_1 & 0
\end{pmatrix}.
\end{equation*}
It is natural to take $Z_{i}$ to be the matrix corresponding to $\zq_{i}{=}1,$ $\zq_{j}{=}0$ for $j{\neq}i$
to define a basis of $T^{1,0}_{\op}M_{0}.$ In an analogous way we can build a basis 
for a linear complement of $\C{H}_{\op}M_{0}$ in $\mathfrak{m}$ by using the $\wq_{i}$'s and $t_{j}$'s. 
\par 
The \textit{vector valued} Levi form $\Lq^{{1,0}}$
is obtained by computing
\begin{equation*}
 \iq\, [Z,\bar{Z}']
\end{equation*}
for matrices $Z,Z'$ of the form $(**)$ 
and e.g. the  $\wq_{1}$ in the entry $(1,2)$ of $(*)$ means that $\Lq^{1,0}(Z_{1},Z_{2})$ is 
a vector in 
$\mathfrak{m}/(\mathfrak{m}\cap(\qt+\bar{\qt}))$  proportional to the projection of 
\begin{equation*} 
 W_{1} = 
\begin{pmatrix}
 0 & 0 & 0 & 0 & 0 & 0\\
 0 & 0 & 0 & 0 & 0 & 0\\
{1} & 0 & 0 & 0 & 0 & 0\\
 0 & 0 & 0 & 0 & 0 & 0\\
 0 & 0 & 0 & 0 & 0 & 0\\
 0 & 0 & 0 & 0 & 0 & 0 
\end{pmatrix}
\end{equation*}
\end{exam}
The entries $\wq_{i}$'s and $t_{j}$'s of $(*)$ can also be taken as dual variables in
$H^{0}_{\op}M.$ In this way $(*)$ can be thought of as the \textit{scalar} Levi form.
We note that it is identically zero on the $3$-plane $\{\wq_{1}=0,\;\wq_{2}=0,t_{2}=0\}$ 
and that all nondegenerate scalar Levi forms have at least one positive and one negative
eigenvalue. 
\theoremstyle{definition}
\newtheorem{exam2}[prop]{Example}
\begin{exam2}
 We consider the flag $M$ of $\SO_7(\C),$ consisting of the projective lines contained in the
 quadric 
\begin{equation*}
 Q=\{[z]\in\CP^6\mid z_1z_7{+}z_2z_6{+}z_3z_5{+}z_4^2{=}0\}
\end{equation*}
 of $\CP^6.$ Take the matrix  
\begin{equation*}
 B= 
\begin{pmatrix}
 0 & 0 & 1\\
 0 & \Id_5 & 0\\
 1 & 0 & 0
\end{pmatrix}.
\end{equation*}
Then  
\begin{equation*}
 \SO(1,6)=\{x{\in}\SO(7,\C)\mid x^*Bx=B\}.
\end{equation*}
In this way, the $\SO(1,6)$-orbit $M_0$ through
the line $\op{=}\{[z]\in\CP^6\mid z\in\langle\eq_1,\eq_2\rangle\}{\subset}Q$ 
is the minimal orbit of $\SO(1,6)$ in $M.$ \par
We have 
\begin{equation*}
 \qt{=}\left.\left\{ 
\begin{pmatrix}
 z_{1,1} & z_{1,2} & z_{1,3} & z_{1,4} & z_{1,5} & z_{1,6} & 0\\
 z_{2,1} & z_{2,2} & z_{2,3} & z_{2,4} & z_{2,5} & 0 & z_{1,7}\\
  0 & 0 & z_{3,3} & z_{3,4} & 0 & z_{3,6} & z_{3,7}\\
   0 & 0 & z_{4,3} & 0 & z_{4,5} & z_{4,6} & z_{4,7}\\
    0 & 0 & 0 & z_{5,4} & z_{5,5} & z_{5,6} & z_{5,7}\\
     0 & 0 & 0 & 0 & 0 & z_{6,6} & z_{6,7}\\
     0 & 0 & 0 & 0 & 0 & z_{7,6} & z_{7,7}  
\end{pmatrix}\; \right| z_{i,j}{=}{-}z_{8-j,8-i}\in\C\right\}.
\end{equation*}
Then, to describe $T_{\op}M_0,$ we construct the matrix of $\ot(1,6)$ which is obtained by adding
to the matrices of the linear complement 
\begin{equation*}
 \mathfrak{m}=\left.\left\{ 
\begin{pmatrix}
 0 & 0 & 0 & 0 & 0 & 0 & 0 \\
  0 & 0 & 0 & 0 & 0 & 0 & 0 \\
   z_{3,1} & z_{3,2} & 0 & 0 & 0 & 0 & 0 \\
    z_{4,1} & z_{4,2} & 0 & 0 & 0 & 0 & 0 \\
     z_{5,1} & z_{5,2} & 0 & 0 & 0 & 0 & 0 \\
      z_{6,1} & 0 & {-}z_{5,2} & {-}z_{4,2} & {-}z_{3,2} & 0 & 0 \\
       0 & {-}z_{6,1} & {-}z_{5,1} & {-}z_{4,1} & {-}z_{3,1} & 0 & 0 \end{pmatrix}
       \; \right| z_{i,j}\in\C\right\}
\end{equation*}
of $\qt$ their conjugate with respect to the real form $\ot(1,6).$ We label by $z_i$'s the entries 
that fall out of $\mathfrak{m}$ by the conjugation and by $w_i$'s and $t_i$'s the complex and real entries
which remain in $\mathfrak{m}$ after conjugation. We obtain in this way the matrices 
\begin{equation*} 
\begin{pmatrix}
 0 & 0 & 0 & 0 & 0 & 0 & 0 \\
 \bar{z}_4 & 0 & {-}\bar{z}_1 & {-}\bar{z}_2 & {-}\bar{z}_3 & 0 & 0 \\
 w & z_1 & 0 & 0 & 0 & \bar{z}_3 & 0\\
 \iq{t} & z_2 & 0 & 0 & 0 & \bar{z}_2 & 0 \\
 \bar{w} & z_3 & 0 & 0 & 0 & \bar{z}_1 & 0 \\
 z_4 & 0 & {-} z_3 & {-} z_2 & {-} z_1 & 0 & 0\\
 0 & {-}z_4 & {-}\bar{w} & {-}\iq{t} & {-}w & {-}\bar{z}_1 & 0
\end{pmatrix}.
\end{equation*}
Then $M_0$ is of type $(4,3)$ and the matrix for its Levi form is 
\begin{equation*} 
\begin{pmatrix}
 0 & 0 & 0 & \wq \\
 0 & 0 & 0 & t\\
 0 & 0 & 0 & \bar{\wq}\\
 \bar{\wq} & t & \wq & 0 
\end{pmatrix},\quad\text{with $\wq\in\C$ and $t{\in}\R.$}
\end{equation*}
In particular, all nonzero scalar Levi forms have signature $+,-,0,0.$
\end{exam2}
As explained in \cite{AMN06}, for homogeneous CR manifolds which,
as the examples discussed above,  
are
orbits of real forms in complex flag manifolds, it is possible to compute
the Levi form by combinatorics on root systems, by using 
the attached cross marked Satake
diagrams. More general examples can be found in \cite{AMN3}.

\end{document}